\documentclass[a4paper]{article}

\usepackage[all]{xy}
\usepackage[latin1]{inputenc}
\usepackage{amscd}
\usepackage{graphicx}
\usepackage{rotating}
\usepackage[english]{babel}
\usepackage{amsmath,amssymb,amsfonts}
\usepackage{setspace}                    
\usepackage{latexsym}
\usepackage{color}
\usepackage{epsfig}

\begin{document}

\title{Some topological  problems on the configuration spaces of Artin and
Coxeter groups\footnote{This paper is an extract of the first part  of the same paper which appeared in the book "Configuration Spaces: Geometry, Combinatorics and Topology", edited by A. Bjorner, F. Cohen, C. DeConcini, C. Procesi, M. Salvetti, Edizioni della Normale, CRM Series, n. 14, 2012 (distributed by Springer). Part of the non-given applications appeared (in a more complete form) in Math. Res. Let., vol. 20, p.1-19, 2013; Atti dell'Accademia Nazionale dei Lincei, vol. 25, p. 233-248, 2014.
}}
\author{D. Moroni \thanks{Institute of Information Science and Technologies
(ISTI), National Research Council of Italy (CNR), Pisa}, M. Salvetti \thanks
{Dipartimento di Matematica, 
Universit\`a di Pisa}, A. Villa \thanks {DIpartimento di Matematica,
Universit\`a di Torino}}
\date {}
\maketitle

\begin{abstract} 
In the first part we review some topological and algebraic aspects in the theory of Artin and Coxeter groups,  both in the finite and infinite case (but still, finitely generated). In the following parts, among other things, we compute the Schwartz genus of the covering associated to the orbit space for all affine Artin groups. We also give a partial computation of the cohomology of the braid group with non-abelian coefficients coming from geometric representations. We introduce an interesting class of "sheaves over posets", which we call "weighted sheaves over posets", and use them for explicit computations.
\end{abstract}

\newcommand{\oC}[0]{\overline{C}}
\newcommand{\oF}[0]{\overline{F}}
\newcommand{\oG}[0]{\overline{G}}
\newcommand{\Y}[0]{\mathbf Y^{(d)}}
\newcommand{\arr}[0]{{\mathcal A}}
\newcommand\X{\mathbf X}
\newcommand\Xd{\mathbf X^{(d)}}
\newcommand\Xdd{\mathbf X^{(d+1)}}
\newcommand\Sd{\mathbf S^{(d)}}
\newcommand\Sdd{\mathbf S^{(d+1)}}
\renewcommand{\S}[0]{\mathbf{S}}
\newcommand\W{\mathbf W}
\newcommand\WG{\mathbf W_{\Gamma}}
\newcommand\WS{\mathbf W^{\Gamma}}
\newcommand\XW{\mathbf X_{\mathbf W}}
\newcommand\Fi{\mathbf \varPhi}
\renewcommand\>{\succ}
\newcommand\po{\lessdot}
\newcommand{\poi}[1]{\lessdot_{#1}^{op}}
\newcommand\Q{\mathbf Q}
\newcommand\Qdd{\mathbf Q^{(d+1)}}
\newcommand\Qd{\mathbf Q^{(d)}}
\newcommand\Qk{\mathbf Q^{(k)}}
\newcommand\join{\bigcup}
\newcommand{\strat}[1]{\Fi}
\newcommand\strk{\Fi ^{(k)}}
\newcommand\strd{\Fi ^{(d)}}
\newcommand\str{\Fi^{(d+1)}}
\newcommand\strY{\Fi^{(d+1)}_0}
\newcommand\strA{\Fi^{(d+1)}_{+}}
\newcommand\strR{\Fi^{(d)}}
\newcommand\prR{pr_{\Re}}
\newcommand\s{\hat}
\newcommand{\Fp}[0]{\Fc^{\prime}}
\newcommand\Fs{\Fc^{\prime\prime}}
\newcommand\cam{\mathbf \varPhi_0}
\newcommand\g{\gamma}
\newcommand\gp{\gamma ^{\prime}}
\newcommand\G{\Gamma}
\newcommand\Gp{\Gamma ^{\prime}}
\newcommand\FO{\{\mathbf 0\}}
\newcommand\GF{ \varGamma}
\newcommand\Gf{\varGamma^{\prime}}
\newcommand\pW{\pi_{\W}}
\newcommand\gk{\hat\g}
\newcommand\GW{\mbf{G_W}}

\newcommand{\A}[0]{\mathcal{A}}
\newcommand{\M}[0]{\mathcal{M}}
\newcommand{\Md}[1]{\mathcal{M}^{(#1)}}
\newcommand{\CS}[2]{\mathcal{#1}^{(#2)}}
\newcommand{\La}[0]{\mathcal{L}(\A)}
\renewcommand{\L}[0]{\mathcal{L}}
\newcommand{\F}[0]{\mathbb{F}}
\newcommand\Fc{\mathcal{F}}
\newcommand\Gc{\mathcal{G}}
\newcommand{\Z}[0]{\mathbb{Z}}
\newcommand{\R}[0]{\mathbb{R}}
\newcommand{\C}[0]{\mathbb{C}}
\newcommand{\CW}[0]{\mathcal{C}}
\newcommand{\N}[0]{\mathbb{N}}
\newcommand{\K}[0]{\mathbb{K}}
\newcommand{\St}[0]{\mathcal S}
\newcommand{\Tbf}[0]{\mathbf{T}}
\newcommand{\Cal}[1]{\mathcal{#1}}
\newcommand{\vs}[0]{\vspace{3mm}}
\newcommand{\ph}[0]{\varphi}
\newcommand{\la}[0]{\lambda}
\newcommand{\bd}[1]{\textbf{#1}}
\newcommand{\e}{\bd e}
\renewcommand{\lessdot}{\vartriangleleft}
\newcommand{\MS}{\tilde{\bd S}}
\newcommand{\te}{\theta}
\newcommand{\Te}{\mathbf{\theta}}
\newcommand{\q}[1]{\mbox{\bfseries{\textit{#1}}}}
\newcommand{\itl}[1]{\textit{#1}}
\newcommand{\p}[0]{p}
\newcommand{\geets}[0]{\longleftarrow}
\newcommand{\too}[0]{\longrightarrow}
\newcommand{\sst}[0]{\subset}
\newcommand{\cl}[1]{\mathcal{#1}}
\newcommand{\into}[0]{\hookrightarrow}
\newcommand{\codim}[0]{\mbox{codim }}
\newcommand{\diag}[0]{\Sigma}
\newcommand{\ea}[0]{\underline{\epsilon}}
\newcommand{\eb}[0]{\overline{\epsilon}}
\newcommand{\im}[0]{\mbox{im }}
\newcommand{\az}[0]{\"{a}}
\newcommand{\ls}[0]{\mathcal L}
\newcommand{\ul}[0]{\underline}
\newcommand{\bin}[2]{  \left(  \begin{array}{c}  #1 \\ #2  \end{array}
  \right)  }
\newcommand{\qbin}[2]{  \left[  \begin{array}{c}  #1 \\ #2
    \end{array}  \right]  }
\newcommand{\qed}[0]{\hspace{\stretch{1}}$\Box$}
\newcommand{\eq}[1][r]
       {\ar@<-3pt>@{-}[#1]
        \ar@<-1pt>@{}[#1]|<{}="gauche"
        \ar@<+0pt>@{}[#1]|-{}="milieu"
        \ar@<+1pt>@{}[#1]|>{}="droite"
        \ar@/^2pt/@{-}"gauche";"milieu"
        \ar@/_2pt/@{-}"milieu";"droite"}
\newcommand{\imm}[1][r] {\ar@{^{(}->}[#1]}
\newcommand{\D}[0]{D^{\scriptscriptstyle{(0)}}}
\newcommand{\dd}[0]{d^{\scriptscriptstyle{(0)}}}
\renewcommand{\ni}[0]{\noindent}
\newcommand{\Ga}[0]{\Gamma}

\newcommand{\flo}[1]{\lfloor #1 \rfloor}
\newcommand{\pr}{\mathrm{Prod}}
\providecommand{\poi}[1]{{\mathbb{W}}_{#1}}
\providecommand{\qed}{{\flushright $\square$\\}}
\renewcommand{\ll}{{ \mathcal{L}_2 (\mathbb{R}^3, d \mu)}}
\providecommand{\p}{{\mathbf{p}}}
\renewcommand{\L}{{\mathbf{L}}}
\renewcommand{\i}{{\mathrm{i}}}
\providecommand{\Br}{{\mathrm{Br}}}
\providecommand{\A}{{\mathbf{A}}}
\newcommand{\lie}[2]{{[ #1 , #2  ]}}
\newcommand{\bra}[2]{{\langle #1, #2 \rangle}}
\newcommand{\no}[1]{{| #1 |}}
\newcommand{\til}[1]{\widetilde{#1}}
\newcommand{\qa}[1]{[#1]}
\newcommand{\fhi}[1]{\{ #1 \}}
\renewcommand{\ln}[1]{\bar{#1}}
\newcommand{\bo}[1]{{\mathcal{#1}}}
\newcommand{\mbf}[1]{{\mathbf{#1}}}

\newcommand{\re}{{\mathbf{R}}}
\newcommand{\B}{{\mathcal{B}}}
\newcommand{\card}[1]{|#1|}
\newcommand{\oth}{\mathrm{otherwise}}
\newcommand{\semidir}{\rtimes}
\newcommand{\de}{\partial}
\newcommand{\di}{{\mathrm{d}}}
\renewcommand{\mod}[1]{\,\mathrm{mod}(#1)}
\newcommand{\compl}{\Xi}
\newcommand{\wI}[0]{\widetilde{I}}
\newcommand{\wqbin}[2]{ \widetilde{\left[ \begin{array}{c} #1 \\ #2
    \end{array} \right] }}
\newcommand{\tss}[0]{\supseteq}
\newcommand{\pmu}[0]{{\pm 1}}
\newcommand{\CA}[0]{\overline{C}}


\newcommand{\br}{\mathbf{G_W}}
\newcommand{\pbr}[1]{PG_{#1}}
\renewcommand{\t}[1]{\tilde{#1}}
\newcommand{\rank}{\mathrm{rk}}
\newcommand{\ind}{\mathrm{Ind}}
\newcommand{\coind}{\mathrm{CoInd}}
\newcommand{\pp}[1]{{(#1)}}


\newtheorem{df}{Definition}[section]
\newtheorem{teo}{Theorem}
\newtheorem{claim}{Claim}
\newtheorem{prop}[df]{Proposition}
\newtheorem{lem}[df]{Lemma}
\newtheorem{cor}[df]{Corollary}
\newtheorem{rmk}[df]{Remark}
\newtheorem{notat}[df]{Notation}
\newtheorem{ex}{Example}
\newtheorem{pf}{Proof}
\newtheorem{for}{}
\newtheorem{conj}{Conjecture}

\newenvironment{es}[1][Example.]{\begin{trivlist}
     \item[\hskip \labelsep {\bfseries #1}]}{\end{trivlist}}
\newenvironment{dm}[1][Proof.]{\begin{trivlist} \item[\hskip
    \labelsep {\bfseries #1}]}{\end{trivlist}}
\newenvironment{dmof}[1][Proof]{\begin{trivlist} \item[\hskip
    \labelsep {\bfseries #1}]}{\end{trivlist}}
\newenvironment{os}[1][Remark.]{\begin{trivlist}
     \item[\hskip \labelsep {\bfseries #1}]}{\end{trivlist}}
\newenvironment{grafi}

\section{Introduction} 
This paper contains some results regarding some particular topological aspects
in the theory of  Artin and Coxeter groups.  Among other things, we compute  the
Schwartz genus of the covering associated to the orbit space for all affine
Artin groups. This result generalizes \cite{salvdec2}, \cite{salvdecproc3}.

Our paper contains also  a  brief review  of some results, concerning the
topology of Artin and Coxeter groups, both in the finite and infinite case (but
still, finitely generated) which are essential to our computations.  
Other reviews, even if considering some interesting aspects of the theory, are
not very satisfactory about the topological underlying structure. Our review will
still be very partial: the more than thirty years old literature on the subject
would require a much longer paper. We concentrate essentially on a single line
of research, which (in our opinion) gives the possibility to produce a 
neat picture of some basic topological situation underlying all the theory in
very few pages.   Such picture is based essentially on \cite{salvettiArtin} (and
\cite{salvDecArtin}) and \cite{salvdec2} (both papers are based on the
construction \cite{salvetti}). Some of the several computations which use these
constructions are cited below.  

For some literature containing many works related to Artin groups the
reader can see the paper  \cite{paris} contained in this book. It should be
added that  almost all people who worked in the theory of Hyperplane
Arrangements have given some contributions to the theory.

The new results that we present here concern some computations of:

\noindent - the cohomology of the braid groups with non-abelian coefficients,
coming from geometric representations of the braid groups into the homology of
an orientable surfaces (see part 3.1);  

\noindent  -the computation of the Schwartz genus of the covering associated to
the orbit space of an affine Artin group (see part 3.2).

We skip the details of the computations of the first application, while we
give  details for the second one.  We introduce here a particularly  interesting
class of ``sheaves over posets'', which we call ``weighted sheaves over
posets''.  We use them for some explicit computations for the top cohomology of
the affine group of type $\t{A}_n$ (the whole cohomology was completely computed
with different methods in \cite{calmorsal1}, \cite{calmorsal2}).  
A natural spectral sequence is associated to such sheaves. We are going to
exploit this general construction in future works.

Some of the ideas which we use here were partially explained in \cite{salvettipreprint}, 
\cite{moronitesi}.
The first of these papers essentially contains a talk given by the second author at a 
Conference in Tokyo, during the 60th-birthday fest in honor of F. Cohen.

\section{General pictures}\label{sec2}
  We will consider finitely generated   {\it Coxeter systems} $(\W,S)$ \ ($S$
finite), so
\begin{equation}\label{presentationW}
\W\ =\ <s\in S\ |\ (ss')^{m(s,s')}=1>
\end{equation}
where $m(s,s')\in\N\cup\{\infty\},\ m(s,s')=m(s',s),\ m(s,s)=1$ 
(for general reference see \cite{bourbaki}, \cite{humphreys}).

\subsection{Case $\W$ finite}
\begin{figure}
\begin{center}
\begin{picture}(400,250)
\includegraphics[width=370pt,height=310pt]{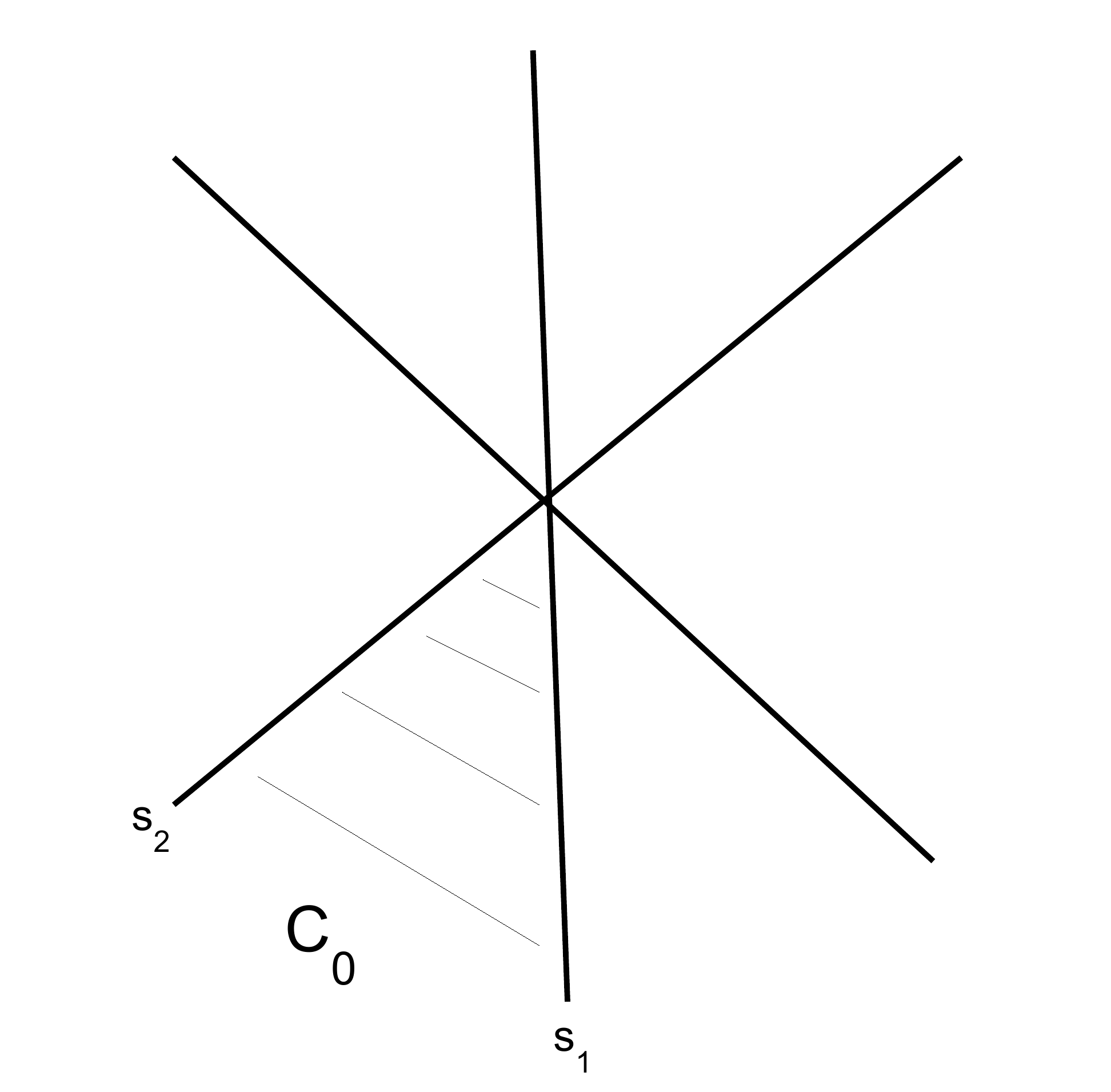}
\end{picture}
\end{center}
\caption{Case $W$ finite} 
\label{fig1} 
\end{figure}

The group $\W$ can be realized as a group generated by reflections in $\R^n$
($n=|S|$).
Let $\A$ be the {\it reflection arrangement}, i.e.
$$\A\ =\ \{H\subset\R^n\ |\ H\ \mbox {is fixed by some reflection in}\ \W\} .$$
Consider also the stratification into {\it facets} $\strat\ := \{ F\}$  of
$\R^n$ induced by $\A$. The codimension-$0$ facets are called {\it chambers}.
They are the connected components of the complement to the arrangement.
All the chambers are simplicial cones, the group acting transitively over the
set of all them.  The Coxeter generator set $S$ corresponds to the set of
reflections with respect to the {\it walls} of a fixed base-chamber $C_0$ (see
fig. \ref{fig1})

Let \ $H_{\C}\ :=\ H+iH\ \subset\ \C^n$ \ be the {\it complexification} of the
hyperplane $H,$ and set
$$\mbf{Y}:=\ \M(\A):=\ \C^n\setminus \cup_{H\in\A}\ H_{\C}$$
be the complement to the complexified arrangement.   The group $W$ acts freely
over $Y$ so we can form the {\it orbit space}
$$\mbf{Y_W}:=\ \mbf{Y/W}$$
which is a complex manifold, actually an affine manifold, by classical results.

Define 
$$\GW:=\ \pi_1(\mbf{Y}_W)$$
as the {\it Artin group} of type $\mbf{W}$ (somebody calls it the {\it
Artin-Brieskorn} group, somebody else calls it the
{\it Artin-Tits} group: for brevity, we just take the intersection of them)
(see \cite{brieskorn, deligne, paris}).
A presentation of $\GW$ is obtained by that of $\W$ by removing the relations
$s^2=1$
\begin{equation}\label{presentationGW}
\GW\ =\ <g_s,\ s\in S\ |\ g_sg_{s'}g_s\dots = g_{s'}g_{s}g_{s'}\dots> 
\end{equation}
(same number   $m(s,s')$ of factors on each side).

From Deligne's theorem (\cite{deligne})  the orbit space $\mbf{Y_W}$ is a space
of type $k(\pi,1),$ so we have
$$H^*(\GW;L)\ =\ H^*(\mbf{Y_W};\L)$$
where $L$ is any $G_W$-module and $\L$ is the corresponding local system over
$\mbf{Y_W}.$

We also recall (\cite{salvettiArtin}, \cite{salvetti}):

\begin{teo}\label{teo:finite}
The orbit space $\mbf{Y_W}$ contracts over a cell complex $\XW$ obtained from a
convex polyhedron $Q$ by explicit identifications over its faces.
\end{teo}

More precisely, one takes one point $x_0\in C_0,$ and set
$$Q:=\  \mbox{convex hull}\ [W.x_0]$$
a polyhedron obtained by the convex hull of the orbit of one point $x_0$ in the
base-chamber. 

For example, in the case of the {\it braid arrangement} one obtains the so
called {\it permutohedron}.

One verifies the following facts.

The $k$-faces of $Q$ are also polyhedra, each of them corresponding to a coset
of a {\it parabolic} subgroup $\WG,$ where $\Ga\subset S$ is a $k$-subset of
$S.$ The correspondence
$$\{ \mbox{faces of $Q$}\} \ \leftrightarrow\ \{ w.\WG,\ \Ga\subset S\}$$
is obtained by taking the  polyhedron given by the convex-hull of the orbit
$\WG.x_0$ and translating it by 
$w.$

One has also (\cite{bourbaki, humphreys}):

\begin{prop}\label{prop:minimal} Inside each coset $w.\WG$ there exists an
unique element  of {\rm minimal length}.
\end{prop}

Here the length is the minimal number of letters (coming from $S$) in a reduced
expression.

For every face $e$ of $Q,$ which corresponds to a coset $w.\WG,$ let
$\beta(e)\in w.\WG$ be the element of minimal length.  Notice that $\W$ permutes
faces of the same dimension. 
Then each pair of faces $e,\ e'$ belonging to the same orbit is identified by
using the homeomorphism
$\beta(e)\beta(e')^{-1}.$

\begin{figure}[h]
\begin{center}
\begin{picture}(370,290)
\includegraphics[width=360pt,height=300pt]{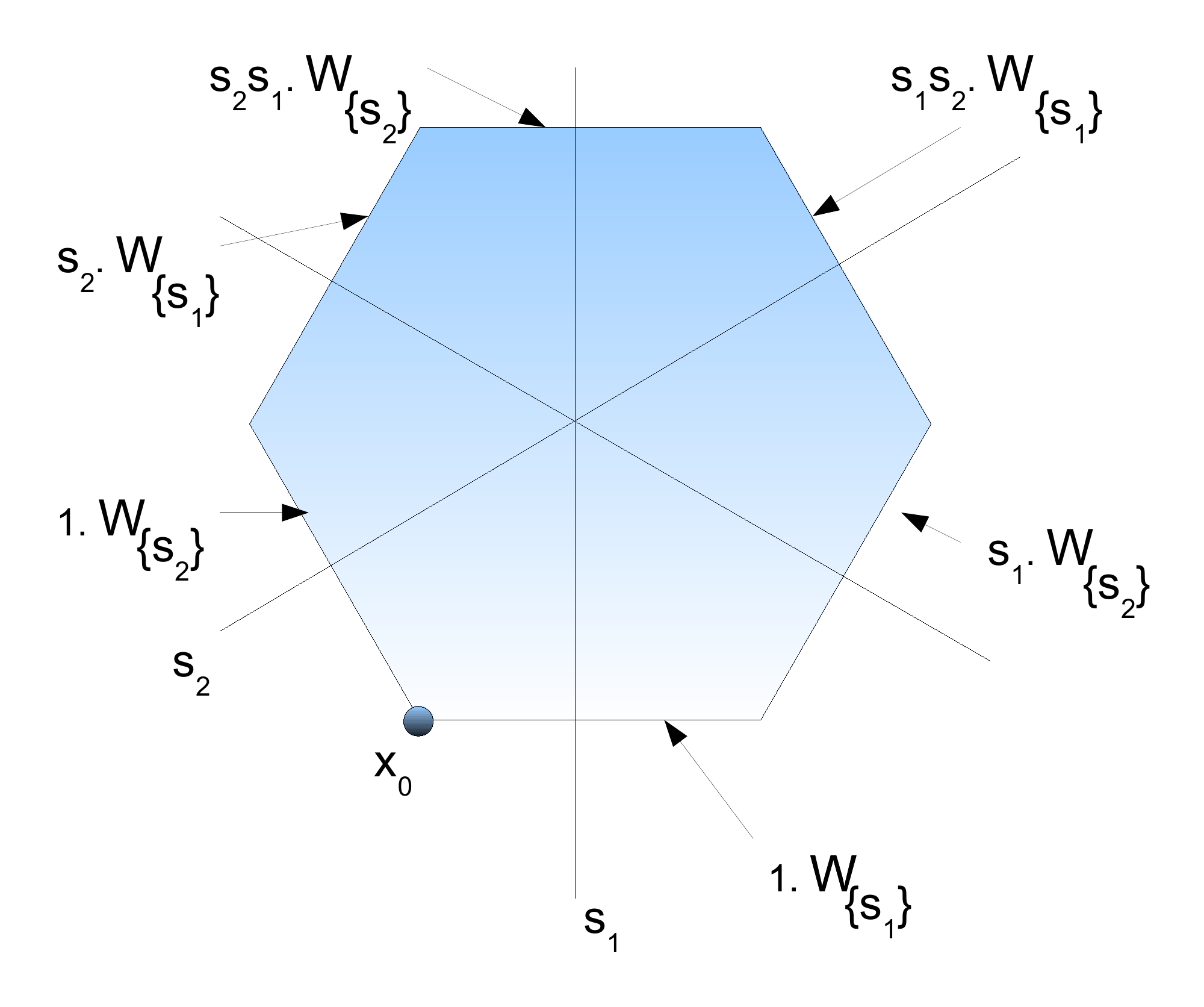}
\end{picture}
\end{center}
\caption{for each $1$-cell $e$ we indicate the corresponding  coset
$\beta(e).\WG$ } 
\label{fig2} 
\end{figure}

\begin{figure}[t]
\begin{center}
\begin{picture}(280,150)
\includegraphics[width=300pt,height=240pt]{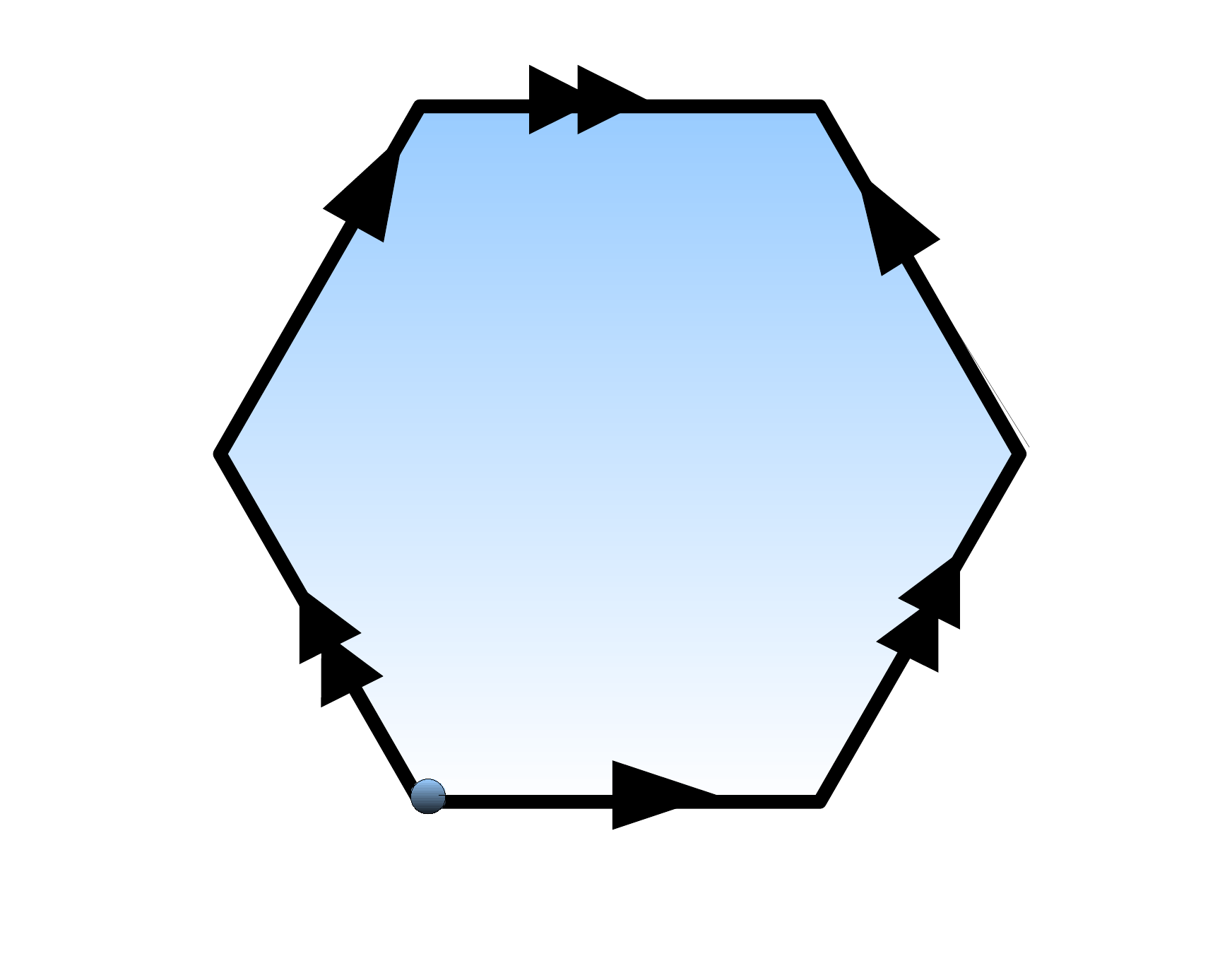}
\end{picture}
\end{center}
\caption{ glue $1$-cells with the same type of arrows } 
\label{fig3} 
\end{figure}

We give in fig. \ref{fig2}, \ref{fig3} the example of the group $\W$ of type
$A_2,$ so $\GW$ is the braid group in three strands. 
The orbit space turns out to have the homotopy type of an hexagon whose edges are
identified according to the given arrows.

\subsection{Case $\W$ infinite}
When $\W$ is infinite (but still $S$ is finite) the theory is analogue with the
following changes
(see \cite{bourbaki,  vinberg}).
One can still realize $\W$ as a group of (non-orthogonal) reflections in $\R^n,$
($n=|S|$) starting from a base chamber $C_0.$ In \cite{bourbaki} the standard
first octant is considered, (so $C_0=\{ x_i>0,\ i=1,\dots,n\}$) but one can
start from a more general 
open cone with vertex $0$ (see \cite{vinberg}). Again, $S$ corresponds to the
set of reflections with respect to the walls of $C_0.$

We recall here the main points which we need.

Let \ $U:= \W.\overline{C}_0$ \ be the orbit of the closure of the base chamber.
$U$ is called the {\it Tits cone} of the Coxeter system. 

Notice that  the closure of the chamber $\overline{C}_0$ is endowed with a
natural stratification into facets (which are still relatively open cones with
vertex $0$). When $C_0$ is the standard positive octant, each facet is given by
imposing some coordinates equal to $0,$ and the remaining coordinates positive.
 
Each reflection in $\W$ is conjugated to a reflection with respect to a wall of
$C_0.$ So,
the arrangement $\A$ of reflection hyperplanes is just the orbit of the walls of
$C_0.$ Each connected component of the complement inside $U$ of the arrangement 
 (again called a chamber) is of
the shape $w.C_0$ for a unique $w\in \W.$   Of course  $\A$  is not locally
finite (e.g. $0$ is contained in all the hyperplanes). The orbits of the facets
of $C_0$ give a ``stratification'' of $U$ into relatively open cells, also
called facets (in general, $U$ is neither open nor closed in $\R^n$).

Recall also (see \cite{bourbaki,vinberg}):
\begin{enumerate}
\item $U$ is a convex cone in $\R^n$ with vertex $0.$
\item $U=\R^n$ \ iff \ $\W$ is finite
\item The stabilizer of a facet $F$ in $U$ is the subgroup $\mbf{W}_F$ generated
by all  the reflections with respect to hyperplanes (in $\A$) containing $F.$ 
So, in general $\W_F$ is not finite. 
\item\label{p3} $U^0:= int(U)$ \ is open in $\R^n$ and a (relatively open) facet
$F\subset\overline{C}_0$ is contained in $U^0$ \ iff \ the stabilizer $\W_F$ is
finite.
\end{enumerate}

 By property \ref{p3} the arrangement  is locally finite in the interior part $U^0.$ 

So we take in this case
$$\mbf{Y}:=\ [U^0\ +\ i\R^n]\ \setminus\ \cup_{H\in\A}\ H_{\C}$$
which corresponds to complexifying only the interior part of the Tits cone.  The
group $\W$ acts (as before) diagonally onto $\mbf{Y},$ and one shows easily
(exactly as in the finite case) that the action is free. 
Therefore, one has an orbit space 
$$\mbf{Y_W}:=\ \mbf{Y/W}$$
which is still a manifold, and a regular covering 
$$\mbf{Y}\ \to\  \mbf{Y_W}$$
with group $\W.$

We still define $\GW\ :=\ \pi_1(\mbf{Y_W})$ as the Artin group of type $\W.$  We
will see in a moment that for  $\GW$ one has a presentation similar to
(\ref{presentationGW}).

In fact, we have very similar constructions. 
\bigskip

Take $x_0\in C_0$ and let $Q$ be the {\it finite} $CW$-complex constructed as
follows.

\noindent For all subsets $\Ga\subset S$  such that $\WG$ is finite, construct a
$|\Ga|$-cell  $Q_{\Ga}$ in $U^0$ as the convex hull of the $\WG$ orbit of $x_0.$
Each $Q_{\Ga}$ is a finite convex polyhedron which contains the point $x_0.$

Let $\X_{\W_{\Ga}}$ be obtained from  $Q_{\Ga}$ by identifications on its faces
defined as in the finite case
(relative to the finite group $\W_{\Ga}$).  Define 
\begin{equation}\label{poly}
Q\ :=\ \cup\ Q_{\Ga}
\end{equation}
(a finite union of convex polyhedra) 
where the union is taken on all the above $\Ga$ for which $\WG$ is finite.
Define also
\begin{equation}\label{salvcomp}
\mbf{X_W}\ :=\ \cup_{\Ga}\ \X_{\W_{\Ga}}.
\end{equation}

\begin{rmk}\label{compatibility}  The definition of $\X_{\W}$ makes sense
because of the following easy fact:

for any common cell $e\subset Q_{\Ga}\cap Q_{\Ga'}$ the minimal element 
$\beta(e)$ is the same when computed in $\WG$ and in $\W_{\Ga'}$.
 \end{rmk}
 
 Moreover, we have the following generalization of the finite case (see
\cite{salvettiArtin})
 
 \begin{teo}\label{teo:main} The $CW$-complex $\X_{\W}$ is deformation retract
of the orbit space $\mbf{Y_W}.$
 \end{teo}

\begin{dm}

Notice that $\mbf{X_W}$ is a finite complex in all cases.
 
First, there exists a regular $CW$-complex $\X\subset \mbf{Y}$ which is
deformation retract of $\mbf{Y},$ and $\X$ is constructed as in \cite{salvetti}.
That paper already worked for the affine cases; however, the same procedure
works in general because we reduce to the locally finite case around faces with finite
stabilizer. 
 
The construction of $\X$ can be chosen invariantly with respect to the action of
$\W,$ which permutes cells of the same dimension.  
 
The action on $\X$ being free,  we look at the orbit space $\mbf{X/W}.$ By
remark  \ref{compatibility}, this reduces to finite cases. \qed\end{dm}

\bigskip

(A similar proof can be obtained by generalizing "combinatorial stratifications" to this situation, 
see \cite{bjorner_ziegler}).

\begin{figure}[h]
\begin{center}
\begin{picture}(340,340)
\includegraphics[width=330pt,height=340pt]{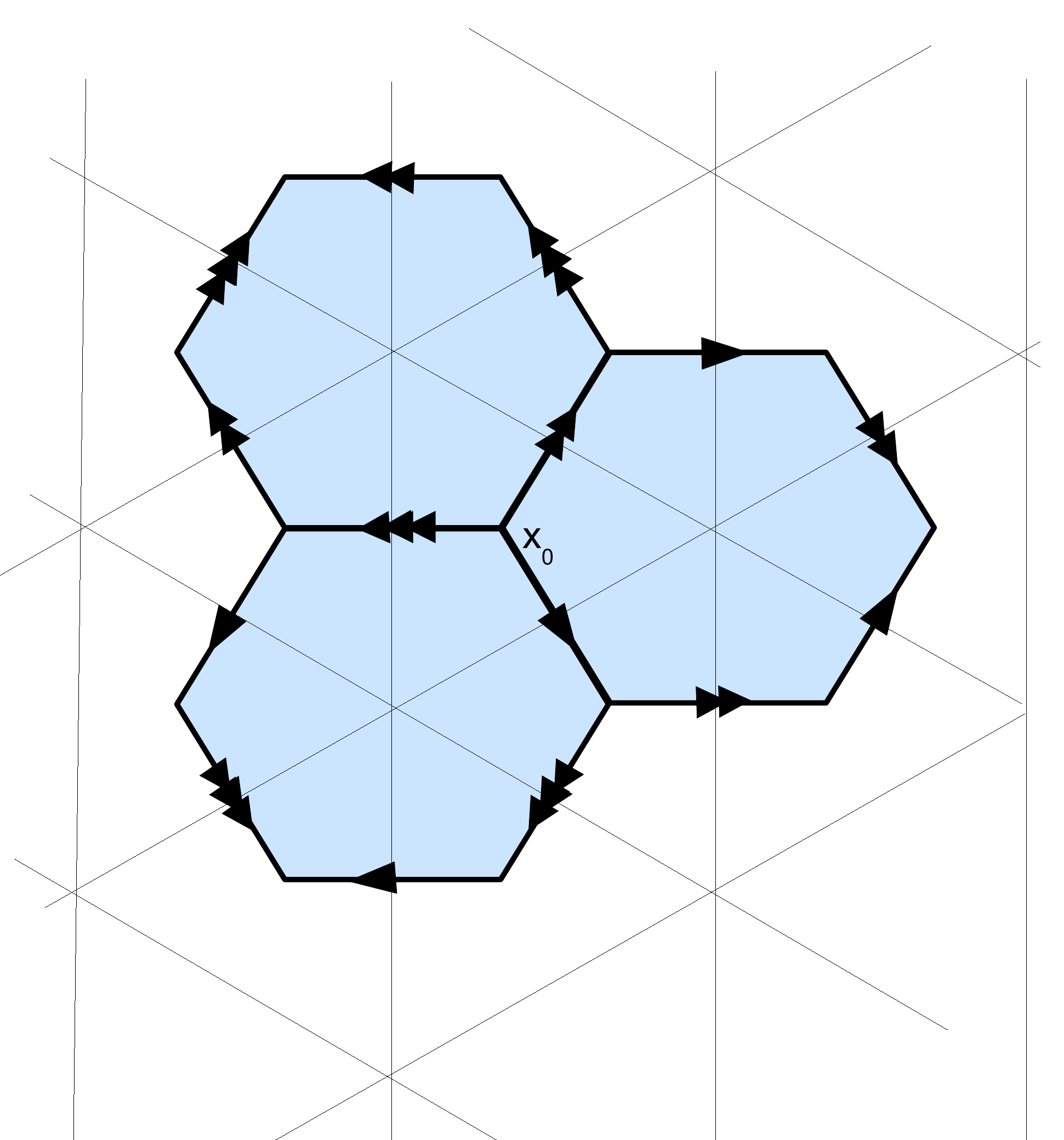}
\end{picture}
\end{center}
\caption{ $Q$ is the union of three hexagons. $\X_{\W}$ is obtained by gluing
$1$-cells with same type of arrow. } 
\label{fig4} 
\end{figure}

Below we give a picture for the case $\tilde{A}_2$  (fig \ref{fig4}).

\begin{conj} The space $\mathbf Y$ is a $k(\pi,1)$ space.
\end{conj}
 
 Of course, if one of the spaces $\mbf{Y}, \ \mbf{Y_W},$ $\X, \ \mbf{X_W}$ is a
$k(\pi,1),$ so are all the other spaces.
 This conjecture is known, besides the finite case, for some affine groups:
$\tilde{A}_n,$ $\tilde{C}_n$ (\cite{Okonek}; in case $\tilde{A}_n$ a different
proof is given in \cite{charney}); case $\tilde{B}_n$ was solved in 
 (\cite{calmorsal3}).
For other known cases see \cite{hendriks}, \cite{charney_davis}.

As an immediate corollary to Theorem \ref{teo:main} we re-find in a very short
way a presentation of the fundamental group
(see \cite{vanderlek}).

\begin{teo}\label{teo:pi1}
A presentation for the Artin group $\GW$ is similar to \ref{presentationGW},
where we have to
consider only pairs $s,\ s'$ such that $m(s,s')$ is finite.
\end{teo}
\begin{dm} We have just to look at the $2$-skeleton of $\X_{\W},$ which is given
as follows. For each pair $s,\ s'$ such that $m(s,s')$ is finite, one has a
$2m(s,s')-$gon with two orbits of edges, which  glue in a similar way as fig.
\ref{fig3}. This gives the relation 
$$g_sg_{s'}g_s\dots\ =\ g_{s'}g_{s}g_{s'}\dots\ \quad \mbox{(} m(s,s')\mbox{ \ 
factors)}.$$  
\qed\end{dm}

\subsection{Algebraic complexes for Artin groups}\label{subsec:Artin}
We refer here mainly to \cite{salvettiArtin},\cite{salvDecArtin}. 

We consider the algebraic complex related to the cell structure of $\X_{\W}.$ It
is given in the following way.

Let $\Z[\GW]$ be the group algebra of $G_W.$  Let $(C_*,\partial_*)$ be the
algebraic complex of free $\Z[\GW]$-modules such that in degree $k$ it is free
with basis $e_J$ corresponding to subsets $J\subset S$ such that $\W_J$ is
finite:

\begin{equation}\label{algcompl}
C_k\ :=\ \oplus_{\begin{array}{l} J\subset S \\ |J|=k \\ \W_J \mbox{ finite} 
\end{array}}\ \Z[\GW]\ e_J
\end{equation}

Let

\begin{equation}\label{bordoArtin}
\partial(e_J)\  :=\ \sum_{\begin{array}{c} I\subset J \\ |I|=k-1\end{array}}\
[I:J] \ T^J_I\ . e_I
\end{equation}
where $[I:J]$ is the incidence number ($=0,\ 1$ or $-1$) of the cells in
$\X_{\W}$ and
$$T^J_I\ :=\ \sum_{\beta\in W^J_I}\ (-1)^{\ell(\beta)}\ g_{\beta}$$
where
\begin{enumerate}
\item $\W^J_I\ :=\ \{\beta\in \W_J\ :\ \ell(\beta \ s)\ >\ \ell (\beta),\
\forall s\ \in\ \W_I \}$  \ 
 is the set of elements of minimal length for  the cosets $\frac{\W_J}{\W_I}$ \
(prop. \ref{prop:minimal});
\item\label{length} if $\beta \in \W^J_I$ and $\beta=s_{i_1}\dots s_{i_k}$ is a
reduced expression then $\ell(\beta)=k;$ 
\item  if $\beta$ is as in \ref{length}  then  $g_{\beta}:= g_{s_{i_1}}\dots
g_{s_{i_k}}.$  One shows that this map 
\begin{equation}\label{section}
\psi: \W \to \GW
\end{equation}
is a well-defined section (not a homomorphism) of the standard surjection 
$\GW \to \W.$ 
\end{enumerate}

\begin{rmk}\label{kappapi1} The orbit space $\X_{\W}$ is a $k(\pi,1)$ iff
$(C_*,\partial_*)$ is acyclic (in this case the augmentation gives a resolution
of $\Z$ into free $\Z[\GW]$-modules).
\end{rmk}
In papers \cite{salvettiArtin}\cite{salvDecArtin} we used the knowledge of the
fact that $\X_{\W}$ is a $k(\pi,1)$ in the finite case to deduce that 
(\ref{algcompl}), (\ref{bordoArtin}) is acyclic. One could try the converse:
prove algebraically that the algebraic complex is exact and conclude that
$\X_{\W}$ is a $k(\pi,1)$  \footnote{Actually this is known in the finite case: an algebraic complex equivalent to the given one was found in \cite{squier}, who proved algebraically that the complex is exact in the finite case.  This paper appeared slightly later than \cite{salvettiArtin}, but actually it was known (it seems, only to the author, who never published something similar) much before because it was essentially the content of his unpublished PhD thesis discussed in 1980. Paper \cite{squier} was written posthumous from his colleagues, who published  the main part of his thesis with some changes (and adding some remarks and theorems which are not contained in the thesis. This whole thing seems very mysterious:  it is not clear why nothing  was published before '94 about Squier's thesis work, as Squier continued to work and did other publications in the 80's. Also, somebody would probably like to know who are the actual authors of the paper appeared in '94.)}.
It is interesting to consider the following abelian representation
(\cite{salvetti}\cite{salvetti_stumbo},
\cite{salvettideconcinistumbo:top}).

Let $R:=A[q,q^{-1}]$ be the ring of Laurent polynomials over a ring $A.$ One can
represent $\GW$ by
\begin{equation}\label{rappr}
g_s\ \mapsto \  [\mbox{multiplication by } -q]  \quad ,\forall s\in S
\end{equation}
($\in Aut(R)$).  This representation, which coincides with the
determinant of the Burau representation, has a very interesting meaning in the
finite case: the cohomology of $\GW$ with coefficient in this representation
equals the trivial cohomology of the Milnor fibre of the associated discriminant
bundle. In other words, the orbit space has a Milnor fibration over $S^1,$ with
Milnor fibre $F,$ and the above twisted cohomology of $\GW$ gives the trivial
cohomology of $F$ as a module over $R,$ the $(-q)$-multiplication corresponding
to the monodromy of the bundle.

The tensor product $C_*\otimes R$ has boundary
\begin{equation}\label{bordoq}
\partial(e_J)\ =\ \sum_{\begin{array}{c}I\subset J\\ |I|=|J|-1\end{array}}\ [I:J]\ \frac{\W_J(q)}{\W_I(q)}\ e_I
\end{equation}
where
$$\W_J(q):=\ \sum_{w\in W_J}\ q^{\textsf{l}(w)}$$
is the Poincar\'e series of the group $\W_J$ (here, a polynomial since the
stabilizers are finite). The denominator $\mathbf W_I(q)$ divides the numerator $\mathbf W_J(q)$
so the quotient is still a polynomial.
 
\begin{es} {\it Case $\mathbf A_n.$}  We have Dynkin graph 

\begin{center}
\begin{tabular}{ccccccccc}
$\circ$ & ----- & $\circ$ & ----- & \dots &  ----- & $\circ$ & ----- &
$\circ$\\ 
1 & & 2 & & & & n-1&  & n
\end{tabular}
\end{center}
Let 
\begin{equation}\label{qanalog}
[k] :=\ \frac{q^k-1}{q-1}\ ;\quad [k]! :=\ \prod_{i=1}^k\ [i]\ ;\ 
\left[\begin{array}{c} k \\ h \end{array} \right] :=\ \frac{[k]!}{[h]![k-h]!}
\end{equation}
For $J\subset S\cong \{ 1,\dots, n \}$ one has
$$\W_J(q)\ =\ \prod_{i=1}^m\ [|\Ga_i(J)|+1]\ !$$
where $\Ga_1(J),\ \Ga_2(J), \dots $ are the connected components of the subgraph
of $\mathbf A_n$ generated by $J$ (fig \ref{fig5}). 

\begin{figure}[h]
\begin{center}
\begin{picture}(450,150)
\includegraphics[width=330pt,height=110pt]{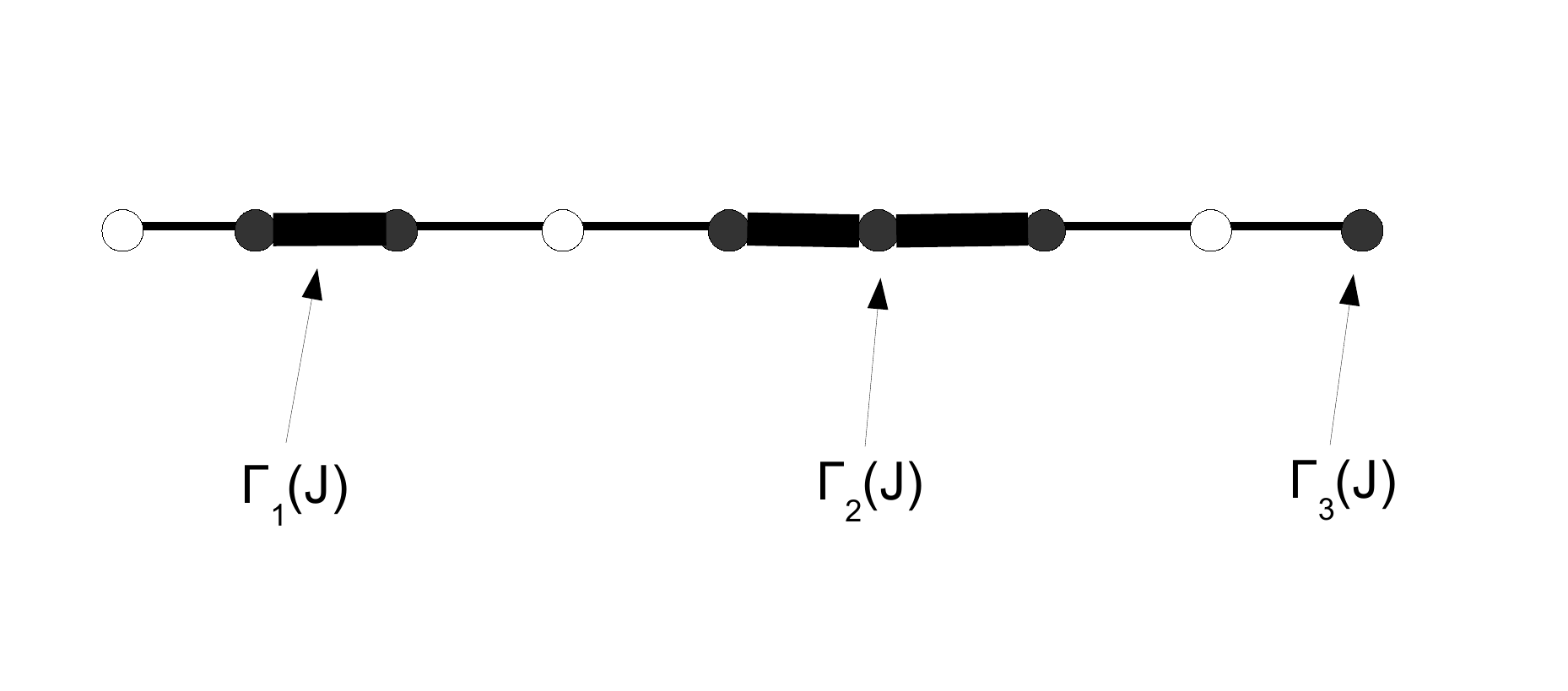}
\end{picture}
\end{center}
\caption{ here $\W_J(q)=[3]![4]![2]!$} 
\label{fig5} 
\end{figure} 
So, all coefficients are of the shape 
$\left[\begin{array}{c} k \\ h \end{array} \right] $ \ ($k,\ h$ depending on
$J,\ I$).
\end{es}

For some computations using these methods for the finite case:

for $H^*(\GW,\Q[q,q^{-1}])$ \cite{arithmetic_braid},  \cite{frenkel} (case
$A_n$)  \cite{arithmetic_artin} (all other cases);
 
for the top cohomology of $H^*(\GW,\Z[q,q^{-1}])$  in all cases
\cite{salvettideconcinistumbo:top}; 

for the cohomology  $H^*(\GW,\Z[q,q^{-1}])$  for all exceptional cases see
\cite{callegarosalvetti};
the same for case $A_n$: see  \cite{callegaro:braid};

In the affine cases: see \cite{calmorsal1} (cohomology in case $\tilde{A}_n$),
\cite{calmorsal2} ($k(\pi,1)$ problem and cohomology for $\tilde{B}_n$).

\subsection{$CW$-complexes for Coxeter groups} 
We refer here essentially to \cite{salvdec2}. 

If $(\W,S)$ is a finite Coxeter system, which is realized as a reflection group
in $\R^n,$  with reflection arrangement  $\A,$ consider the subspace arrangement
in $\R^{nd}\ \cong (\R^n)^d$   given by
$$\A^{(d)}:=\ \{ H^{(d)}\ \}$$
where $H^{(d)}$ is the codimensional-$d$ subspace given by
``$d$-complexification'' of the hyperplane $H\in \A:$
$$H^{(d)}:=\ \{ (X_1,\dots,X_d)\ :\ X_i\in\R^n,\ X_i\in\ H\ \}.$$
For $d=2$ we have the standard complexification of the hyperplanes.
Let 
$$\Y:=\ \R^{nd}\setminus \cup_{H\in\A}\ H^{(d)}.$$
As before, the group $W$ acts freely on $Y^{(d)}$ and we consider the orbit
space
$$\Y_{\W}\ :=\ \Y/\W.$$
Recall:
\begin{teo}\label{kpi1} The space
$$\mbf{Y}_{\mbf{W}}^{(\infty)}:=\ \left[ \lim_{\stackrel{\longrightarrow}{d}}\ \Y\right]
/ W\ =\  \left[ \lim_{\stackrel{\longrightarrow}{d}}\ \Y_{\W}\right] $$
is a space of type $k(W,1).$
\qed\end{teo}

In case $(W,S)$ is infinite (still, $S$ finite) one has to substitute 
$\R^n\times\dots\times\R^n$ ($d$ factors) with the space
$$\mathcal U_W^{(d)}:= U_0\times\R^n\times\dots\times\R^n$$
($d-1$ factors equal to $\R^n$), and the space $Y^{(d)}$ becomes
\begin{equation}\label{confspace}
\Y:= \ \mathcal U_W^{(d)} \setminus \cup_{H\in\A}\ H^{(d)}.
\end{equation}

\begin{teo}\label{kpi1bis} For finitely generated $\W,$ the same conclusion as
in theorem \ref{kpi1} holds, by  
taking definition (\ref{confspace}) for $Y^{(d)}.$
\end{teo}
So, different from the case of Artin groups, we always get a $k(\pi,1)$ space
here.

Recall also the construction of a $CW$-complex which generalizes that given for
Artin groups. 
\begin{teo}\label{S-complex}  When $\W$ is finite, the space $\Y_{\W}$ contracts
over a $CW$-complex $\X^{(d)}_{\W}$
such that
$$\{ k-\mbox{cells of } \X^{(d)}_{\W}\}\ \longleftrightarrow \{ \mbox{flags }
\mbf{\Gamma}:=(\Ga_1\supset\dots \supset\Ga_d):\ \Ga_1\subset S, \ \sum_{i= 1}^d
|\Ga_i|\ =\ k\} $$
Passing to the limit, $\mbf{Y}^{(\infty)}_{\W}=k(W,1)$ contracts over a
$CW$-complex $\X^{(\infty)}_{\W}$
such that
$$\{ k-\mbox{cells of } \X^{(\infty)}_{\W}\}\ \longleftrightarrow \{ \mbox{flags
} \mbf{\Gamma}:=(\Ga_1\supset\Ga_2\supset\dots):\ \Ga_1\subset S, \ \sum_{i\geq
1} |\Ga_i|\ =\ k\} $$
\qed\end{teo}

Notice that $\X^{(\infty)}_{\W}$ does not have finite dimension but the number
of  $k$-cells is finite, given by $\left({n+k-1}\atop{k}\right).$ 

The case $\W$ infinite has to be modified as in case of Artin groups, by
considering flags composed with subsets $\Ga\subset S$ such that $\W_{\Ga}$ is
finite. In the limit, the theorem is

\begin{teo}\label{S-complexbis}  For $S$ finite, the space 
$\mbf{Y}^{(\infty)}_{\W}=k(W,1)$ contracts over a $CW$-complex
$\X^{(\infty)}_{\W}$ such that the set of $k-cells$ of  $\X^{(\infty)}_{\W}\}$
corresponds to
$$\{ \mbox{flags } \mbf{\Gamma}:=(\Ga_1\supset\Ga_2\supset\dots):\ \Ga_1\subset
S,\ \W_{\Ga_1}\mbox{ finite, }\ \sum_{i\geq 1} |\Ga_i|\ =\ k\} $$
\end{teo}

\subsection{Algebraic complexes for Coxeter groups}\label{Coxeter}
Consider the algebraic complex
$(C_*^\pp{d}, \de)$ of free  
$\Z[\W]$-modules, where
$$
C_k^\pp{d}:= \bigoplus_{\scriptsize \begin{array}{c} \mbf{\Gamma}:
\sum_{1}^d \card{\Gamma_i}=k \\ \card{W_{\Gamma_1}}< \infty
\end{array}} \Z[\W] e(\mbf{\Gamma})
$$
The generators of $C_*$ are in one to
one correspondence with the cells of $\X^{(\infty)}_{\W}$. The
expression of the boundary is the following:
$$
\de e(\mbf{\Gamma})= \sum_{ \scriptsize
                            \begin{array}{c}
                                         1 \leq i \leq d \\
                                    \card{\Gamma_i}>\card{\Gamma_{i+1}}
                            \end{array}
                            }
                    \sum_{\tau \in \Gamma_i}
                    \sum_{ \scriptsize
                            \begin{array}{c}
                                    \beta \in \W^{\Gamma_i \setminus \{ \tau
                                    \}}_{\Gamma_i}\\
                                    \beta^{-1} \Gamma_{i+1} \beta
                                    \sst \Gamma_i \setminus \{ \tau
                                    \}
                            \end{array}
                            }
                            (-1)^{\alpha(\mbf{\Gamma}, i, \tau,
                            \beta)}\beta e(\mbf{\Gamma'})
$$
where
$$
\mbf{\Gamma'}= (\Gamma_1\supset \ldots\supset \Gamma_{i-1}\supset \Gamma_i \setminus
\{\tau\} \supset \beta^{-1} \Gamma_{i+1} \beta\supset \ldots \supset
\beta^{-1}\Gamma_d \beta )
$$
and $(-1)^{\alpha(\mbf{\Gamma}, i, \tau, \beta)}$ is an incidence
index. To get a precise expression for ${\alpha(\mbf{\Gamma},i,
\tau, \beta)}$, fix a linear order on $S$ and let
\begin{align*}
\mu(\Gamma_i, \tau):=&\card{j \in \Gamma \textrm{ s.t. } j \leq \tau}\\
\sigma(\beta, \Gamma_j):=& \card{(a,b) \in \Gamma_j\times \Gamma_j
\textrm{ s.t. } a<b \textrm{ and } \beta(a)>\beta(b)}
\end{align*}
in other words, $\mu(\Gamma_i, \tau)$ is the number of reflections
in $\Gamma_i$ less or equal to $\tau$ and $\sigma(\beta,
\Gamma_j)$ is the number of inversions operated by $\beta$ on
$\Gamma_j$. Then we define:
$$
    {\alpha(\mbf{\Gamma}, i, \tau, \beta)}= i \ell(\beta) +
    \sum_{j=1}^{i-1} \card{\Gamma_j} + \mu(\Gamma_i, \tau) +
    \sum_{j=i+1}^{d} \sigma(\beta, \Gamma_j)
$$
where $\ell$ is the length function in the Coxeter group.

Let now  \ $C_*=:\lim_{d \rightarrow \infty} C_*^\pp{d}$. The flags in $C_k$ are
the (infinite) sequences 
$$\mbf{\Gamma}= (\Gamma_1\supset \Gamma_{2}\supset \dots ) $$
such that  $\sum_{i\geq 1} |\Ga_i|\ =\ k,$ \ so they still have a finite number
of nonempty $\Ga_i.$

\begin{teo}\label{algcompCox}
For any finitely generated $\W,$ the algebraic complex $(C_*,\de_*)$  gives a
free resolution of the trivial 
$\Z[\W]$-module $\Z.$
\end{teo}
The proof follows straightforward from the remark that the limit space
$\mbf{Y}^{(\infty)}_{\mbf W},$  so 
$\mbf{X}^{(\infty)}_{\mbf W},$
is a space of type $k(\pi,1).$

\section{Applications} \label{sec3}
\bigskip

\centerline{....... \emph{continue in the paper cited in the footnote in the first page} ......................}
\bigskip

\bibliographystyle{amsalpha}
\bibliography{bibliomat}

\providecommand{\bysame}{\leavevmode\hbox to3em{\hrulefill}\thinspace}
\providecommand{\MR}{\relax\ifhmode\unskip\space\fi MR }
\providecommand{\MRhref}[2]{%
  \href{http://www.ams.org/mathscinet-getitem?mr=#1}{#2}
}
\providecommand{\href}[2]{#2}
\begin{thebibliography}{CMS08b}

\bibitem[Bou68]{bourbaki}
N.~Bourbaki, \emph{Groupes et algebr\`es de {L}ie}, vol. Chapters {IV-VI},
  Hermann, 1968.

\bibitem[Bri71]{brieskorn}
E.~Brieskorn, \emph{Die {F}undamentalgruppe des {R}aumes der regul{\"a}ren
  {O}rbits einer endlichen komplexen {S}piegelunsgruppe}, Invent. Math.
  \textbf{12} (1971), 57--61.

\bibitem[BZ92]{bjorner_ziegler}
A.~{Bjorner} and G.~Ziegler, \emph{Combinatorial stratifications of complex
  arrangements}, Jour. Amer. Math. Soc. \textbf{5} (1992), 105--149.

\bibitem[Cal06]{callegaro:braid}
F.~Callegaro, \emph{The homology of the {M}ilnor fiber for classical braid
  groups}, Alg. Geom. Top. \textbf{6} (2006), 1903--1923.

\bibitem[CD95]{charney_davis}
R.~Charney and M.W. Davis, \emph{The $k(\pi,1)$-problem for hyperplane
  complements associated to infinite reflection groups}, J. of AMS \textbf{8}
  (1995), 597--627.

\bibitem[Cha95]{charney}
R.~Charney, \emph{Geodesic automation and growth functions for {A}rtin groups
  of finite type}, Math. Ann. \textbf{301} (1995), 307--324.

\bibitem[CMS08a]{calmorsal2}
F.~Callegaro, D.~Moroni, and M.~Salvetti, \emph{Cohomology of affine {A}rtin
  groups and applications}, Trans. Amer. Mat. Soc. \textbf{360} (2008),
  4169--4188.

\bibitem[CMS08b]{calmorsal1}
\bysame, \emph{Cohomology of {A}rtin groups of type $\tilde{A}_n, \tilde{B}_n$
  and applications}, Geom. \& Top. Mon. \textbf{13} (2008), 85--104.

\bibitem[CMS10]{calmorsal3}
\bysame, \emph{The ${K}(\pi,1)$ problem for the affine {A}rtin group of type
  $\tilde{B}_n$ and its cohomology,}, Jour. Eur. Math. Soc. \textbf{12} (2010),
  1--22.

\bibitem[CS04]{callegarosalvetti}
F.~Callegaro and M.~Salvetti, \emph{Integral cohomology of the {M}ilnor fibre
  of the discriminant bundle associated with a finite {C}oxeter group}, C. R.
  Acad. Sci. Paris, Ser. I \textbf{339} (2004), 573--578.

\bibitem[Del72]{deligne}
P.~Deligne, \emph{Les immeubles des groupes de tresses g\'en\'eralis\'es},
  Inventiones math. \textbf{17} (1972), 273--302.

\bibitem[DPS01]{arithmetic_braid}
C.~{De Concini}, C.~Procesi, and M.~Salvetti, \emph{Arithmetic properties of
  the cohomology of braid groups}, Topology \textbf{40} (2001), no.~4,
  739--751.

\bibitem[DPS04]{salvdecproc3}
\bysame, \emph{On the equation of degree $6$}, Comm. Math. Helv. \textbf{79}
  (2004), 605--617.

\bibitem[DPSS99]{arithmetic_artin}
C.~{De Concini}, C.~Procesi, M.~Salvetti, and F.~Stumbo, \emph{Arithmetic
  properties of the cohomology of {A}rtin groups}, Ann. Scuola Norm. Sup. Pisa
  Cl. Sci. \textbf{28} (1999), no.~4, 695--717.

\bibitem[DS96]{salvDecArtin}
C.~{De Concini} and M.~Salvetti, \emph{Cohomology of {A}rtin groups}, Math.
  Res. Lett. \textbf{3} (1996), 293--297.

\bibitem[DS00]{salvdec2}
\bysame, \emph{Cohomology of {A}rtin groups and {C}oxeter groups}, Math. Res.
  Lett. \textbf{7} (2000), 213--232.

\bibitem[DSS97]{salvettideconcinistumbo:top}
C.~{De Concini}, M.~Salvetti, and F.~Stumbo, \emph{The top-cohomology of
  {Artin} groups with coefficients in rank-1 local systems over $\mathbb{Z}$},
  Topology and its Applications \textbf{78} (1997), 5--20.

\bibitem[Fre88]{frenkel}
E.~V. Frenkel, \emph{Cohomology of the commutator subgroup of the braids
  group}, Func. Anal. Appl. \textbf{22} (1988), no.~3, 248--250.

\bibitem[Hen85]{hendriks}
H.~Hendriks, \emph{Hyperplane complements of large type}, Invent. Math.
  \textbf{79} (1985), 375--381.

\bibitem[Hum90]{humphreys}
J.E. Humphreys, \emph{Reflection groups and {C}oxeter groups}, Cambridge
  University Press, 1990.

\bibitem[L.P11]{paris}
L.Paris, \emph{Basic question on \mbox{Artin-Tits} groups}, these Proceedings
  (2011).

\bibitem[Mor06]{moronitesi}
D.~Moroni, \emph{Finite and infinite type artin groups: Topological aspects and
  cohomological computations}, PhD thesis (2006).

\bibitem[Oko79]{Okonek}
C.~Okonek, \emph{Das $k(\pi, 1)$-problem f{\"ur} die affinen wurzelsysteme vom
  typ {$A_n$, $C_n$}}, Mathematische Zeitschrift \textbf{168} (1979), 143--148.

\bibitem[Sal87]{salvetti}
M.~Salvetti, \emph{Topology of the complement of real hyperplanes in
  $\mathbb{C}^n$}, Invent. Math. \textbf{88} (1987), no.~3, 603--618.

\bibitem[Sal94]{salvettiArtin}
\bysame, \emph{The homotopy type of {A}rtin groups}, Math. Res. Lett.
  \textbf{1} (1994), 567--577.

\bibitem[Sal05]{salvettipreprint}
\bysame, \emph{On the {C}ohomology and {T}opology of {A}rtin and {C}oxeter
  groups}, Pubblicazioni Dipartimento di Matematica L.Tonelli, Pisa (2005).

\bibitem[Squ94]{squier}
C.C. Squier, \emph{The homological algebra of {A}rtin groups}, Math. Scand.
  \textbf{75} (1994), no.~1, 5--43.

\bibitem[SS97]{salvetti_stumbo}
M.~Salvetti and F.~Stumbo, \emph{Artin groups associated to infinite {C}oxeter
  groups}, Discrete Mathematics \textbf{163} (1997), 129--138.

\bibitem[vdL83]{vanderlek}
H.~van~der Lek, \emph{The homotopy type of complex hyperplane complements},
  Ph.D. thesis, University of Nijmegan, 1983.

\bibitem[Vin71]{vinberg}
E.B. Vinberg, \emph{Discrete linear groups generated by reflections}, Math.
  USSR Izvestija \textbf{5} (1971), no.~5, 1083--1119.

\end{thebibliography}

\end{document}